\documentclass{article}
\usepackage{multicol,graphicx,color}
\usepackage{pslatex}
\usepackage{amsthm}
\usepackage{amsmath}
\usepackage{pstricks}

\newcommand{\myfig}[3][0]{
\begin{center}
  \vspace{0.2cm}
  \includegraphics[width=#3\hsize,angle=#1]{#2}

  \nobreak\medskip
\end{center}
}

\newcommand{\mycaption}[1]{
  \vspace{0.2cm}
  \begin{quote}
    {{\sc Figure} \arabic{figure}: #1}
  \end{quote}
  \vspace{0.2cm}
  \stepcounter{figure}
}

\theoremstyle{plain}
\newtheorem{theorem}{Theorem}
\newtheorem{corollary}[theorem]{Corollary}

\author{Ouyang Tiancheng \ , \ Duokui Yan\\
        Department of Mathematics, Brigham Young University\\
}

\title{Periodic Solutions with Alternating Singularities in the Collinear Four-body Problem}

\begin{document}

\maketitle

\begin{abstract}
This paper shows the existence of a periodic orbit with singularity in the symmetric collinear four body problem. In each period of the orbit, there is a binary collision (BC) between the inner two bodies and a simultaneous binary collision (SBC) of the two clusters on both sides of the origin.
The system is regularized and the existence is proven by using the implicit function theorem and a continuity argument on differential equations of the regularized Hamiltonian.
\end{abstract}
 \textbf{ Key Words:} N-body Problem, Collinear Motion, Binary Collision, Simultaneous Binary Collision, Periodic Solution with Singularity

\section{Introduction}
In this paper, we prove the existence of a periodic orbit with singularities in the symmetric collinear four body problem of celestial mechanics.
The orbit regularly alternates between two kinds of collisions: a binary collision between the inner two bodies and a simultaneous
binary collision of the two clusters on both sides of the origin. Many numerical results related to similar periodic orbits have been published.
This type of orbit was first discovered in the collinear three-body problem with equal mass by Schubart \cite{SC} (1956) and subsequently extended to unequal masses by (H\'{e}non \cite{HE2}, 1977)
who had also studied the continuation of the family into planar motions \cite{HE} (H\'{e}non, 1976). The unequal mass case was further
 investigated numerically by Mikkola and Hietarinta \cite{HM2} (1991). There is one such interplay orbit for each mass ratio of the collinear
 three-body problem, which when stable is surrounded by a local region of phase space filled with quasiperiodic orbits.
 These are one of the dominant features of the dynamics as illustrated by Hietarinta and Mikkola \cite{HM} (1993) .\\

 Recently, two mathematical proofs of the existence of Schubart orbits have been provided: one is a topological proof given by Mockel\cite{MO}, 2008 and the other is a variational proof given by Venturelli\cite{VE}, 2008. In our paper, a new mathematical proof is given for the Schubart-like orbit in the symmetric collinear four body problem. \\
 Inspired by Sweatman's work\cite{SW} and \cite{SW2}, we consider the following symmetric case: the four masses are $1$, $m$, $m$, and $1$ respectively, and they are symmetric about the center of the mass, which is assumed to be the origin. By analyzing the regularized
Hamiltonian, we use the implicit theorem and a continuity argument to show the existence of this type of periodic orbit. In
section 2.1, we give an important estimate for the maximum
distance of the outer bodies such that the motion of the inner
bodies is monotonic. In section 2.2, we introduce the regularized
Hamiltonian. In section 3, the implicit function theorem is used to
show that the net momentum of the cluster at simultaneous binary collision (SBC) is continuous with
respect to one particular variable $R$, which is the square root of
the maximum distance of the outer bodies. We show that there exist
distinct values of $R$ such that the net momentum of the cluster at
SBC is positive for one value of $R$ and negative for the other.
Therefore, by continuity there exists some $R$ such that the net
momentum of the cluster at SBC is 0, which guarantees the existence
of the periodic orbit because of symmetry.
\section{The Setting and the orbit}
\subsection{The setting in Cartesian Coordinate System}
\begin{center}
\psset{xunit=1in,yunit=1in}
\begin{pspicture}(0,-0.4)(2,.6)

 \psline{->}(-1,0)(3,0)
\psdots[dotsize=3pt 1](-.8, 0)(0, 0)(1.5,0)(2.3,0)
\psdots[dotsize=5pt 1](.75, 0)
 \rput(-.8, -0.15){$-x_1$}
 \rput(0, -0.15){$-x_2$}
 \rput(0.75, -0.15){$0$}
 \rput(1.5, -0.15){$x_2$}
 \rput(2.3, -0.15){$x_1$}
  \rput(-.8, .15){$1$}
 \rput(0, .15){$m$}

 \rput(1.5, 0.15){$m$}
 \rput(2.3, 0.15){$1$}
 \rput(1,-0.4){Fig. 1}
\end{pspicture}
\end{center}
From right to left, number the four bodies from 1 to 4. As in the
problem setting, the masses for the four bodies are $1, \ m, \ m,
\text{and} \ 1$ respectively. The system remains symmetrically
distributed about the center of mass.  The coordinates for the
bodies are $x_1$, $x_2$, $-x_2$ and $-x_1$ respectively. And body 1
to 4 have velocities $\dot{x}_1$,
$\dot{x}_2$,$-\dot{x}_2$,$-\dot{x}_1$ respectively.

The Netowanian equations are
\begin{equation*}
\ddot{x}_1= -\frac{1}{4 x_1^2} -\frac{m}{(x_1+x_2)^2} -
\frac{m}{(x_1-x_2)^2}
\end{equation*}

\begin{equation*}
\ddot{x}_2= -\frac{m}{4 x_2^2} -\frac{1}{(x_1+x_2)^2} +
\frac{1}{(x_1-x_2)^2}
\end{equation*}

We will choose the total energy $E=-1$.\\

In this paper, we are interested in finding a special periodic orbit with singularities.
The orbit alternates between binary collision (BC) between the inner two bodies 2 and 3, and SBC between bodies 1 and 2 and bodies 3 and 4. By introducing a new set of transformations, the singularities of BC and SBC in our problem can be regularized.

\subsection{The Setting in the Regularized System}\label{0}
We will adopt Sweatman's\cite{SW2} work to regularize the system.The system has Hamiltonian
$$H= \frac{1}{4} w_1^2+ \frac{1}{4m} w_2^2 - \frac{1}{2x_1} - \frac{m^2}{2x_2} - \frac{2m}{x_1+x_2}-\frac{2m}{x_1-x_2},  $$
where  $w_1= 2 \dot{x}_1$ and $w_2= 2 m \dot{x}_2$ are the conjugate momenta to $x_1$ and $x_2$. Introduce a canonical transformation
$$q_1=x_1-x_2,  \ \ \ \ \ \ \  q_2=2 x_2, \ \ \ \ \ \ \  p_1=w_1,  \ \ \ \ \ \ \  p_2= \frac{1}{2}(w_1+w_2).$$
This results in a new form for the Hamiltonian
$$H= (1+\frac{1}{m})\frac{p_1^2}{4}   -\frac{p_1p_2}{m} +\frac{p_2^2}{m} -\frac{2m}{q_1}- \frac{m^2}{q_2} - \frac{2m}{q_1+q_2}- \frac{1}{2q_1+q_2}. $$

To regularize the equations of motion, we introduce a Levi-Civita type of  canonical transformation
$$ Q _i= \sqrt{q_i} ,  \ \ \ \ \ \ \  P_i= 2 Q_i p_i \ \ (i=1,2),$$
and we also replace time $t$ by the new independent variable $s$ which is given by $\frac{d t}{d s}= q_1q_2$. In the extended phase space, this produces a regularized Hamiltonian
\begin{align*}
 \Gamma & = \frac{d t}{d s}  (H-E) \\
& = \frac{1}{16} Q_2^2 P_1^2 + \frac{Q_2^2P_1^2- 4Q_1 Q_2 P_1 P_2+ 4Q_1^2P_2^2}{16m}-m^2Q_1^2-2m Q_2^2 -\frac{2m Q_1^2 Q_2^2}{Q_1^2+ Q_2^2}-\frac{ Q_1^2 Q_2^2}{2Q_1^2+ Q_2^2} - Q_1^2 Q_2^2 E,
\end{align*}
where $E=-1$ is the total energy.\\

We start at BC with initial conditions
$$x_1(0)=A, \qquad  x_2(0)=0, \qquad  \dot{x}_1(0)=0, \qquad \dot{x}_2(0)= +\infty$$
which is a singular point.\\

To analyze the motion, it is necessary to deal with the singularity in the regularized coordinate system. The corresponding initial conditions at $s=0$ in this new coordinate system are:
$$Q_1(0)=R, \qquad Q_2(0)=0,  \qquad P_1(0)=0,  \qquad P_2(0)=2 m^{\frac{3}{2}} ,$$
where $A=R^2$.\\

\subsection{Estimation of A}
Intuitively, if $A$ is big enough, there will be multiple BCs before the first SBC happens. In order to find the desired orbit, we will have to give an estimation of $A$ such that there will be no BC for  $t\in (0, t_1]$, where $t_1$ is the time of first SBC.\\
\textbf{Definition:} Let a body has velocity $v=0$ at time $t^{*}$. If there exists an interval $[t_a, t_b]$, such that $t_a<t^{*}<t_b$, and $v$ is positive for $t \in [t_a, t^{*})$ and is negative for  $t \in (t^{*}, t_b]$, or $v$ is negative for $t \in [t_a, t^{*})$ and is positive for  $t \in (t^{*}, t_b]$, then we call $t^{*}$ is the turning time and the position of the body at $t^{*}$ is called the turning point.

\begin{center}
\psset{xunit=1in,yunit=1in}
\begin{pspicture}(0,-0.4)(2,.6)

 \psline{->}(-1,0)(3,0)
\psdots[dotsize=3pt 1](-.8, 0)(0, 0)(1.5,0)(2.3,0)
\psdots[dotsize=5pt 1](.75, 0)
 \rput(-.8, -0.15){$-x_1$}
 \rput(0, -0.15){$-x_2$}
 \rput(0.75, -0.15){$0$}
 \rput(1.5, -0.15){$x_2$}
 \rput(2.3, -0.15){$x_1$}
\psarc[linewidth=1pt](1.45, 0.2){0.118}{270}{450}
\psline[linewidth=1.1pt](0.78, 0.15)(1.45, 0.15)
\psline[linewidth=1.1pt]{<-}(1.3, .25)(1.45, .25)
\psline[linewidth=0.3pt]{<-}(1.52, 0.23)(1.7, 0.4)
\rput(1.86, 0.4){\scriptsize{Turning Point}}

\psarc[linewidth=1pt](0.05, 0.2){0.116}{90}{270}
\psline[linewidth=1.1pt](0.72, 0.15)(0.05, 0.15)
\psline[linewidth=1.1pt]{->}(0.05, .25)(0.2, .25)

\psline[linewidth=1pt]{->}(-0.8, .15)(-0.6, .15)
\psline[linewidth=1pt]{->}(2.3, .15)(2.1, .15)

\end{pspicture}
\end{center}

\begin{theorem}\label{11}
There exists an $A_0$, such that the second body has no turning point for $t \in(0, t_1]$ whenever $0<A \leq A_0$, where $t_1$ is
the time when the first SBC happens. Further, the second body will have at least one turning point for $t \in(0, t_1]$ if $ A > A_0$.
\end{theorem}

\begin{proof}
Consider the Newtonian equations of $x_1$ and $x_2$:
\begin{equation}\label{x1}
\ddot{x}_1= -\frac{1}{4 x_1^2} -\frac{m}{(x_1+x_2)^2} -
\frac{m}{(x_1-x_2)^2}
\end{equation}
\begin{equation}\label{x2}
\ddot{x}_2= -\frac{m}{4 x_2^2} -\frac{1}{(x_1+x_2)^2} +
\frac{1}{(x_1-x_2)^2}
\end{equation}

In order to get an upper bound of $A$, we consider a necessary for
the
second body changing direction. Assume $t=t^{*}<t_1$ is the time when  $\dot{x}_2=0$, then $\ddot{x}_2 \leq 0$ for $t\in (0, t^{*}]$.\\

 Let $x_1(t^{*})=ax_2(t^{*})$ with $a>1$. Because $\ddot{x}_2(t^{*}) \leq 0$,  $a$ must satisfy

 \begin{equation*}
-\frac{m}{4} -\frac{1}{(a+1)^2} + \frac{1}{(a-1)^2} \leq 0
 \end{equation*}
i.e.
  \begin{equation}\label{a}
  16a\leq m(a^2-1)^2.
  \end{equation}

 Also, for $t\in [0, t^{*}]$, $x1/x2 \geq x_1(t^{*})/x_2(t^{*})=a$. By the setting of the problem, $a>1$.\\

Rewrite equation \ref{x1} as:

$$ -\ddot{x}_1= \frac{1}{4 x_1^2} +\frac{m}{(x_1+x_2)^2} + \frac{m}{(x_1-x_2)^2} $$
$$ = \frac{1}{x_1^2} \left[ \frac{1}{4} + \frac{2m(1+\frac{x_2^2}{x_1^2})}{(1-\frac{x_2^2}{x_1^2} )^2 } \right].$$

Note $2m(1+y)/(1-y )^2 $ is increasing with respect to $y$
if $y<1$. In our case, $y= x_2^2/x_1^2  \leq
1/a^2<1$, then

$$ \frac{1}{x_1^2} \left[ \frac{1}{4} + \frac{2m(1+\frac{x_2^2}{x_1^2})}{(1-\frac{x_2^2}{x_1^2} )^2 } \right] \leq\frac{1}{x_1^2} \left[ \frac{1}{4} + \frac{2(1+\frac{1}{a^2})}{(1-\frac{1}{a^2} )^2 } \right]  $$
$$=\frac{1}{x_1^2} \left[\frac{1}{4} + \frac{2ma^2(a^2+1)}{(a^2-1 )^2 } \right] $$                                                                                                                           Therefore,
$$  -\ddot{x}_1 \leq \frac{1}{x_1^2} \left[ \frac{1}{4} + \frac{2a^2(a^2+1)}{(a^2-1 )^2 } \right] .$$

Let $x_1(0)=A$ and $x_1(t^{*})=A_1<A$. Since $x_1$ is
decreasing for $t \in [0, t^{*}]$, i.e. $\dot{x}_1 < 0$  as $t \in [0, t^{*}]$.
The following inequality is true for $t \in [0, t^{*}]$:
$$
-\dot{x}_1 \frac{1}{x_1^2} \left[ \frac{1}{4} +
\frac{2ma^2(a^2+1)}{(a^2-1 )^2 } \right]  \geq \dot{x}_1
\ddot{x}_1,$$
integrate this from $t=0$ to $t=t^{*}$ to get:
$$\left[ \frac{1}{4} +
\frac{2ma^2(a^2+1)}{(a^2-1 )^2 } \right] \int_0^{t^{*}}  \left(-\dot{x}_1
\frac{1}{x_1^2}\right) dt \geq  \int_0^{t^{*}} \dot{x}_1 \ddot{x}_1
dt,$$
$$(\frac{1}{A_1}- \frac{1}{A}) \left[
\frac{1}{4} + \frac{2 m a^2(a^2+1)}{(a^2-1 )^2 } \right]\geq
\frac{1}{2}\dot{x}_1^2 (t^{*}). $$ Then
\begin{equation}\label{b}
(\frac{1}{A_1}- \frac{1}{A}) \left[ \frac{1}{2} + \frac{4 m
a^2(a^2+1)}{(a^2-1 )^2 } \right]\geq \dot{x}_1^2 (t^{*}).
\end{equation}
As $E=-1$, then
$$-1= \dot{x}_1^2 (t^{*})  - \left[\frac{1}{2x_1} +\frac{m^2}{2x_2} + \frac{2m}{x_1+x_2} + \frac{2m}{x_1-x_2}   \right] $$
$$= \dot{x}_1^2 (t^{*}) -\frac{1}{A_1} \left[ \frac{1}{2} + \frac{m^2 a}{2} +\frac{2 m}{1+ \frac{1}{a}} + \frac{2 m}{1- \frac{1}{a}}\right] $$
$$= \dot{x}_1^2 (t^{*}) -\frac{1}{A_1} \left[ \frac{1}{2} + \frac{m^2 a}{2} +\frac{4 m a^2}{a^2-1}  \right] .$$
By inequality \ref{b},
$$-1 \leq (\frac{1}{A_1}- \frac{1}{A}) \left[ \frac{1}{2} + \frac{4 m a^2(a^2+1)}{(a^2-1 )^2 } \right] - \frac{1}{A_1} \left[ \frac{1}{2} + \frac{m^2 a}{2} +\frac{4m a^2}{a^2-1}  \right] $$
  $$= \frac{1}{A_1} \left[ \frac{4 m a^2(a^2+1)}{(a^2-1 )^2 }-  \frac{m^2 a}{2} -\frac{4m a^2}{a^2-1} \right] - \frac{1}{A} \left[ \frac{1}{2} + \frac{4m a^2(a^2+1)}{(a^2-1 )^2 } \right]$$

 $$=\frac{1}{A_1} \frac{m a[16a-m(a^2-1)^2] }{2(a^2-1)^2} - \frac{1}{A} \left[ \frac{1}{2} + \frac{4 m a^2(a^2+1)}{(a^2-1 )^2 } \right] . $$

 Applying inequality \ref{a} to this gives,
 $$-1\leq  \frac{1}{A_1} \frac{m a[16a-m (a^2-1)^2] }{2(a^2-1)^2} - \frac{1}{A} \left[\frac{1}{2} + \frac{4 ma^2(a^2+1)}{(a^2-1 )^2 } \right] \leq - \frac{1}{A} \left[ \frac{1}{2} + \frac{4 ma^2(a^2+1)}{(a^2-1 )^2 } \right] .$$
  Then
\begin{equation}\label{c}
  A \geq \frac{1}{2} + \frac{4 m a^2(a^2+1)}{(a^2-1 )^2 }.
\end{equation}

By inequalities \ref{a} and \ref{c},
$$A > \frac{1}{2}+ 4m + \frac{12 m}{ a^2-1}+ \frac{8 m}{ (a^2-1)^2} , $$
where  $a^4-2a^2-\frac{16a}{m} +1 \geq 0$. In other words, if the following two inequalities $A\leq 1/2+ 4m + 12 m/(a^2-1)+ 8 m/ (a^2-1)^2 $ and $a^4-2a^2-16a/m +1 \leq 0$ hold, $x_2$ will be monotonically increasing as $t\in (0, t_1]$. Choose $A_0$, such that $\dot{x}_2(t^{*})= \ddot{x}_2 (t^{*})=0 $ for some time $t^{*} \in (0, t_1]$, that is
\begin{equation}\label{4}
a^4-2a^2-\frac{16a}{m} +1 = 0
\end{equation}
and
$$A_0 \geq \frac{1}{2}+ 4m + \frac{12 m}{ a^2-1}+ \frac{8 m}{ (a^2-1)^2} .$$
When $A=A_0$, derivative the equation \ref{x2} with respect to $t$ evaluated at $t^{*}$:
$$\dddot{x}_2(t^{*})= \frac{m \dot{x}_2(t^{*})}{2x_2^3(t^{*})} + \frac{2(\dot{x}_1(t^{*})+\dot{x}_2(t^{*}))}{\left[x_1(t^{*})+ x_2(t^{*})\right]^3} -\frac{2(\dot{x}_1(t^{*})-\dot{x}_2(t^{*}))}{\left[x_1(t^{*})- x_2(t^{*})\right]^3} .$$
Since \ $\dot{x}_2(t^{*})=0$, \ $\ddot{x}_2(t^{*})=0$,\
$\dot{x}_1(t^{*})<0$, \  $x_1(t^{*})>x_2(t^{*})>0$,
$$\dddot{x}_2(t^{*})= 2\dot{x}_1(t^{*}) \left[\frac{1}{\left[x_1(t^{*})+ x_2(t^{*})\right]^3}- \frac{1}{\left[x_1(t^{*})- x_2(t^{*})\right]^3}\right]>0.$$
Hence, when $A=A_0$, $x_2$ will increase monotonically for $t\in(0, t_1]$. Therefore, $x_2$ is increasing for $t\in (0,t_1]$ whenever $A \leq A_0$, where $t_1$ is the time of the first SBC, and $x_2$ will have at least one local maximum for $t\in (0,t_0)$ when $A>A_0$.
\end{proof}
\textbf{Remark:} In the special case $m=1$, equation \ref{4} becomes $a^4-2a^2-16a +1 = 0$ with $a>1$, $a \approx 2.766.$ Then
$$A_0 \geq \frac{1}{2}+\frac{4a^2(a^2+1)}{(a^2-1 )^2 }>6.485. $$

\section{Existence of the Periodic Orbit}

Recall that in Section \ref{0},
\begin{align*}
 \Gamma & = \frac{d t}{d s}  (H-E) \\
& = \frac{1}{16} Q_2^2 P_1^2 + \frac{Q_2^2P_1^2- 4Q_1 Q_2 P_1 P_2+ 4Q_1^2P_2^2}{16m}-m^2Q_1^2-2m Q_2^2 -\frac{2m Q_1^2 Q_2^2}{Q_1^2+ Q_2^2}-\frac{ Q_1^2 Q_2^2}{2Q_1^2+ Q_2^2} - Q_1^2 Q_2^2 E,
\end{align*}
where $E=-1$ is the total energy. The initial conditions at $s=0$ are:
$$Q_1(0)=R, \ \ \ \ \ \ \ \ \ Q_2(0)=0,  \ \ \ \ \ \ \ \ \  P_1(0)=0,  \ \ \ \ \ \ \ \ \  P_2(0)=2m^{\frac{3}{2}} .$$
By Theorem \ref{11}, when $0<R= \sqrt{A} \leq \sqrt{A_0} $, $Q_2= \sqrt{2x_2}$ is increasing from $s=0$ to $s=s_1$, where $s_1$ is the time when the first SBC happens.\\

The equations of motion from the regularized Hamiltonian $\Gamma$ are:
\begin{equation}\label{e1}
  Q'_1=  \frac{1+m}{8m} Q_2^2P_1-\frac{1}{4m} Q_1Q_2P_2,
\end{equation}

\begin{equation}\label{e2}
  Q'_2= \frac{1}{2m} Q_1^2P_2- \frac{1}{4m} Q_1 Q_2 P_2 ,
\end{equation}

\begin{equation}\label{e3}
 P'_1= \frac{1}{4m} P_1 P_2 Q_2- \frac{1}{2m}  Q_1 P_2^2  + 2m^2 Q_1 + \frac{4 m Q_1 Q_2^4}{( Q_1^2+ Q_2^2)^2} + \frac{2 Q_1 Q_2^4}{(2 Q_1^2+ Q_2^2)^2} - 2Q_1 Q_2 ^2,
\end{equation}

\begin{equation}\label{e4}
 P'_2=\frac{1}{4m} P_1 P_2 Q_1- \frac{1+m}{8m} Q_2 P_1^2  + 4mQ_2 + \frac{4m Q_1^4 Q_2}{ (Q_1^2+ Q_2^2)^2} + \frac{4 Q_1 ^4 Q_2}{(2 Q_1^2+ Q_2^2)^2} - 2Q_1^2 Q_2.
\end{equation}
where $'$ is the derivative with respect to $s$. The initial conditions are
$$ Q_1(0)= R, \qquad Q_2(0)=0, \qquad  P_1(0)=0, \qquad P_2(0)= 2. $$
At the time $s_1$ when the first SBC happens,
$$Q_1(s_1)=0,  \qquad Q_2(s_1)=R_1, \qquad P_1(s_1)=-\frac{8m}{\sqrt{2m+2}}.$$

To prove the existence of the periodic orbit, we are going to find a value of $R$, such that
$P_2(s_1)=0$ when $Q_1(s_1)=0$.\\

\begin{theorem}\label{1}
For the regularized Hamiltonian $\Gamma$, let the initial
condition be $P_1(0)=0 $, $P_2(0)=2$, $Q_1(0)=R>0$ and $Q_2(0)=0$,
where $R\in(0,   \sqrt{A_0}].$ Assume at time $s_1=s_1(R)$, $Q_2>0$,
$Q_1=0$, which means SBC occurs. Also, assume $Q_2>0$ for $0<s \leq
s_1. $ Then $P_2 (s_1, R)$ is a continuous function of $R$.
\end{theorem}

\begin{proof}
Since the Hamiltonian $\Gamma$ is regularized, the solution $P_i=P_i(s, R)$ and $Q_i=Q_i(s, R)$ are continuous functions with respect to $s$ and $R$. We are going to show $s_1=s_1(R)$ is a continuous function of $R$. In order to apply the implicit function theorem for $Q_1=Q_1(s_1, R)=0$, we need to show that $(\partial Q_1 / \partial s) (s_1, R) \neq 0.$\\
By the regularized Hamiltonian $\Gamma$,
$$\frac{\partial Q_1}{\partial s } \mid_{(s_1, R)} = \Gamma_{P_1} \mid_{(s_1, R)}=  \left[\frac{1+m}{8m} Q_2^2P_1-\frac{1}{4m} Q_1Q_2P_2\right] \mid_{(s_1, R)}  $$

Note that for fixed $R$, $\Gamma=0$ at any time $s$. At $s=s_1$, $Q_1=Q_1(s_1, R)= 0$, then $P_1=P_1(s_1, R)= 4$. Therefore,
$$\frac{\partial Q_1}{\partial s } \mid_{(s_1, R)} = \left[\frac{1+m}{8m} Q_2^2P_1-\frac{1}{4m} Q_1Q_2P_2\right] \mid_{(s_1, R)}= \sqrt{\frac{m+1 }{2}}Q^2 _{2} (s_1, R) >0. $$

By the implicit function theorem, $s_1$ is a continuous function of $R$. Then  $P_2 (s_1, R)$ is also a continuous function of $R$.
\end{proof}

\begin{corollary}
There exists $R$ such that $P_2 (s_1)=P_2 (s_1, R)=0$.
\end{corollary}

\begin{proof}
First, we are going to show that there exists an $R>0$ such that $P_2(s_1)>0$.\\

From equations \ref{e1}-\ref{e4},
$$(P_1 Q_1+ P_2Q_2)'= P'_1 Q_1+ P_1 Q'_1+ P'_2 Q_2+ P_2 Q'_2$$
$$= 4 m Q_2^2+ 2 m^2 Q_1^2 + 2Q_1^2Q_2^2 \left[\frac{2 m }{ Q_1^2+ Q_2^2}+ \frac{1}{ 2Q_1^2+ Q_2^2} -2 \right]$$

Note that $A>x_1(t)>x_2(t)>0$ for $t\in [0, t_1]$, then $0\leq Q_1 \leq R= \sqrt{A}$ and $2x_1(t_1)=2x_2(t_1)= Q_2^2(s_1) =R_1^2 < 2A= 2R^2$, that is $0< Q_2 < R_1< \sqrt{2} R$ for $s\in [0, s_1]$. Thus

$$\frac{2m}{ Q_1^2+ Q_2^2}+ \frac{1}{ 2Q_1^2+ Q_2^2} \geq \frac{2 m}{3R^2} + \frac{1}{4R^2} .$$

Choose $R= \sqrt{\frac{m}{3}}$,

$$\frac{2m }{ Q_1^2+ Q_2^2}+ \frac{1}{ 2Q_1^2+ Q_2^2}-2 \geq \frac{3}{4m} >0 ,$$
so $(P_1 Q_1+ P_2Q_2)' \geq 0$ for $s \in [0, s_1]$, or  $P_1 Q_1+ P_2Q_2$ is increasing  for $s \in [0, s_1]$. \\
From the initial conditions,
$$\left(P_1 Q_1+ P_2 Q_2 \right) \mid _{s=0}= 0 .$$
Hence,
$$0 < \left( P_1 Q_1+ P_2Q_2\right) \mid _{s=s_1}= R_1 P_2(s_1).$$
Therefore, when $R=\sqrt{\frac{m}{3}}$,  \ $ P_2(s_1) > 0$ .\\

Next, we will show $P_2(s_1)$ is negative when $R^2=A_0$.\\

At $A=A_0$, by the proof of theorem \ref{11}, there exists a time $t^{*}<t_1$, such that $\dot{x}_2(t^{*})=0$ and $\dot{x}_1(t^{*})<0$. Then $ \dot{x}_1(t^{*})+\dot{x}_2(t^{*})<0$. Consider the sum of the Newtonian equations \ref{x1} and \ref{x2}:

$$\ddot{x_1}+ \ddot{x_2} = -\frac{1}{4x_1^2} -\frac{1}{4x_2^2} - \frac{2}{(x_1+x_2)^2} <0 ,$$
which means $\dot{x}_1(t)+\dot{x}_2(t)$ is a decreasing function with respect to $t$. Hence, $\dot{x}_1(t_1)+\dot{x}_2(t_1)<\dot{x}_1(t^{*})+\dot{x}_2(t^{*})<0$. Note that $P_2(s_1)/[2Q_2(s_1)]= p_2(t_1)= \dot{x}_1(t_1)+\dot{x}_2(t_1)<0$, and $Q_2(s_1) >0$, so $P_2(s_1)<0$.\\

Therefore, by continuity, there must exist an $R$, such that $P_2(s_1)=0$ where $s_1$ is the time when the first SBC happens.

\end{proof}

\begin{theorem}
If $R$ satisfies $P_2 (s_1, R)=0$, then the orbit will be a one-dimensional Schubart-like periodic orbit.
\end{theorem}

\begin{proof}
At time $s=0$, a BC happens between bodies 2 and 3. At time $s=s_1$, a SBC occurs.
Since the system is regularized, the solution $\{ P_i, Q_i \} \ \ (i=1,2)$  is continuous. \\

At time $s=0$, $$Q_1(0)=R, \ \ \ \ \ Q_2(0)=0, \ \ \ \ \  P_1(0)=0 , \ \ \ \ \ P_2(0)=2 m^{\frac{3}{2}}.$$

At time $s=s_1$,
$$Q_1(s_1)=0, \ \ \ \ \ Q_2(s_1)= R_1, \ \ \ \ \  P_1(s_1)=-\frac{8m}{\sqrt{2m+2}} , \ \ \ \ \ P_2(s_1)=0,$$
where $R_1$ is a positive number. From the Hamiltonian $\Gamma$, we can see that $Q'_1(s_1)=\frac{1}{4} Q_2 (Q_2 P_1 -Q_1 P_2) <0$, $Q'_2(s_1)=\frac{1}{4} Q_1 (2 Q_1 P_2 -Q_2 P_1) =0$.
This means that $Q_2(s_1)$ is the local maximum of $Q_2$. In other words, when time passes $s_1$, $Q_2$ will decrease. Similarly, $Q_1$ will decrease when time passes $s_1$. \\
At the time  $s=s_2$ when the second SBC occurs,
$$Q_1(s_2)=-R, \ \ \ \ \ Q_2(s_2)= 0, \ \ \ \ \  P_1(s_2)=0 , \ \ \ \ \ P_2(s_2)=-2 m^{\frac{3}{2}}.$$
Compare the motion for $s\in  [0, s_1]$ and the motion for $s\in  [s_1, s_2]$. By the uniqueness of the regularized Hamiltonian system and symmetry, the orbit from $s=s_1$ to $s= s_2$ will be the same trajectory from $s=s_1$ to $s=0$ by reversing the direction of each velocity. Further, $s_2=2s_1$. By symmetry and uniqueness again, at time $s=3s_1$,
$$Q_1(3s_1)=0, \ \ \ \ \ Q_2(3s_1)= -R_1, \ \ \ \ \  P_1(3s_1)= \frac{8m}{\sqrt{2m+2}}  , \ \ \ \ \ P_2(3s_1)=0.$$
At time $s=4s_1$,
$$Q_1(4s_1)=R, \ \ \ \ \ Q_2(4s_1)=0, \ \ \ \ \  P_1(4s_1)=0 , \ \ \ \ \ P_2(4s_1)=2 m^{\frac{3}{2}},$$ which is exactly the same as the initial condition at $s=0$. Then the orbit from $s=0$ to $s=4s_1$ will generate one period.

\end{proof}

The following figure is a picture of the periodic solution for $m=1$ in terms of $\{Q_1, Q_2, P_1, P_2\}$. The initial conditions are

$$ Q_1(0)= 2.295, \qquad Q_2(0)=0, \qquad  P_1(0)=0, \qquad P_2(0)= 2. $$

\begin{center}
\myfig{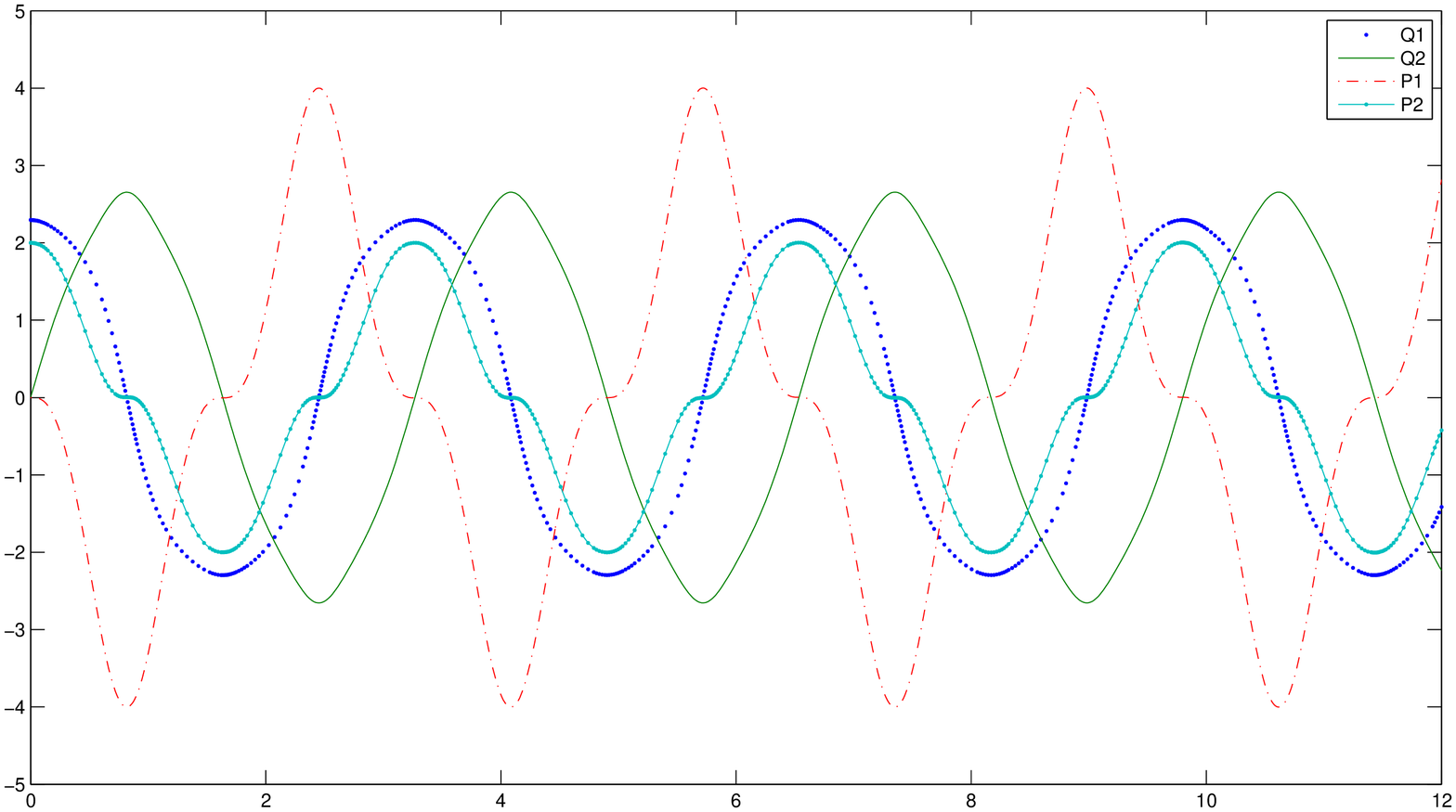}{1 }
\mycaption{The horizontal axis represents time $s$. }
\end{center}


\section{Acknowledgements}

We are pleased to acknowledge several valuable conversations with Dr. Lennard Bakker on these and related topics. We are also indebted to Dr. Lennard Bakker for reading the original manuscript and suggesting improvements.


\end{document}